# ON ARITHMETIC KLEINIAN GROUPS GENERATED BY THREE HALF-TURNS

MICHAEL BELOLIPETSKY

## 1. Introduction.

We study a generalization of the Fuchsian triangle groups to the hyperbolic 3-space $\mathbb{H}^3$. These are groups generated by half-turns in three hyperbolic lines in $\mathbb{H}^3$. If the lines lie in a single plane and intersect each other then the problem reduces to a Fuchsian triangle group generated by reflections of the plane in sides of a hyperbolic triangle. In general case the three lines with their common perpendiculars give rise to a hyperbolic right-angled hexagon. We look at the right-angled hexagons in $\mathbb{H}^3$ as spacial analogs of hyperbolic triangles. This way of thinking has been suggested by W. Fenchel and proved to be fruitful in hyperbolic geometry (see [2]). From technical point of view, our approach to studying generalized triangle groups is close to the methods from H. Helling, J. Mennicke and E.B. Vinberg paper [9] in which another generalization of triangle groups was considered.

The final purpose of our work is to characterize all arithmetic generalized triangle groups. This problem is closely connected with the problem of classification of the arithmetic 2-generator Kleinian groups lately studied by C. Maclachlan and G. Martin [10]. We explain this connection in the second section of this paper.

From the section 3 we begin studying groups generated by three half-turns. We start with finding out an appropriate matrix representation of the generalized triangle group in $\mathrm{SL}(2,\mathbb{C})$. Here, in order to parametrize our class of groups by three complex numbers, we make use of the notion of complex distance between hyperbolic lines and an expression of the distance in terms of line-matrices [2]. In order to prove that the parameters determine the group up to conjugation (Lemma 3.1) we refer to uniqueness theorems for hyperbolic right-angled hexagons [2].

In section 4 we introduce an arithmeticity test for the groups generated by three half-turns. Our test is based on well-known criterions of arithmeticity for two-generated Kleinian groups [8], [4].

In section 5 we turn to a partial case of the regular triangle groups which are generated by three half-turns with equal complex distances between the axes. For this case we prove a finiteness theorem for the number of the arithmetic groups. To eliminate groups which split into free products we make use of geometric estimates supplied by Klein combination theorem. We also describe parameters of non cocompact arithmetic regular triangle groups and present singular sets of corresponding hyperbolic orbifolds.

In the last section we return to a discussion of the general case which still keeps many open problems. Even the finiteness theorem for the number of arithmetic generalized triangle groups, as the connected finiteness result for the number of arithmetic 2-generated Kleinian groups, has not been proved or refuted yet.

Author thanks Alexander D. Mednykh for suggesting the problem and many helpful discussions.

## 2. Two-generator Kleinian groups and groups generated by three half-turns.

We begin with recalling some basic definitions. Following W. Fenchel [2] we consider hyperbolic 3-space as a space of j-quaternions $\mathbb{H}^3 = \{w = z + \zeta j \mid z \in \mathbb{C}, \zeta > 0\}$ endowed with Riemannian metric $ds = |dw|/\zeta$ of constant negative curvature. This way we can associate to a matrix $\pm \begin{pmatrix} a_{11} & a_{12} \\ a_{21} & a_{22} \end{pmatrix} \in \mathrm{PSL}(2,\mathbb{C})$ Möbius transformation $w = z + \zeta j \to (a_{11}w + a_{12})(a_{21}w + a_{22})^{-1}$ and obtain an action of $\mathrm{PSL}(2,\mathbb{C})$ as a full group $Iso^+(\mathbb{H}^3)$ of orientation preserving isometries of hyperbolic space.

**Definition 2.1.** *A discrete nonelementary subgroup of $Iso^+(\mathbb{H}^3)$ is called* Kleinian *group. Here* nonelementary *means that the group does not have a finite index abelian subgroup. The elements of a discrete subgroup of $Iso^+(\mathbb{H}^3)$, other then the identity, fall into three types:*

  parabolic *elements are those which have only one fixed point in $\overline{\mathbb{H}^3} = \mathbb{H}^3 \cup \overline{\mathbb{C}}$;*

  loxodromic *elements have precisely two fixed points in $\overline{\mathbb{H}^3}$. Hyperbolic line joining the points is invariant under the loxodromic isometry and is called its axis;*

  elliptic *elements fix a line in $\mathbb{H}^3$. The fixed line is called an axis of the elliptic isometry.*

Let $\Gamma_0 = \langle L, M \rangle$ be a discrete subgroup of $\mathrm{PSL}(2,\mathbb{C})$. There is a natural way to associate to $\Gamma_0$ a three-generator subgroup $\Gamma$ of $\mathrm{PSL}(2,\mathbb{C})$ which is equal to $\Gamma_0$ or contains it as a subgroup of index 2. From the arithmetic point of view there is no difference between the group and its subgroup or extension in $\mathrm{PSL}(2,\mathbb{C})$ of finite index [8]. So the following construction makes it possible to study arbitrary two-generator arithmetic Kleinian groups using the groups generated by three half-turns in hyperbolic space.

The group $\Gamma$ is constructed as follows (see [5]): If neither $L$ nor $M$ are parabolic then let $N$ be a common perpendicular to the axes of $L$ and $M$. In case of parabolic generator corresponding end of the line $N$ is the fixed point of the parabolic. Obviously, this construction uniquely defines line $N$ for any Kleinian group $\Gamma_0 = \langle L, M \rangle$. Let $c$ be the half-turn about $N$. Then we can factor $L$ as $ac$ and $M$ as $cb$ ([2], p. 47). We obtain the group $\Gamma = \langle a, b, c \rangle$ which contains $\Gamma_0$ as a subgroup of index at most 2.

Visa versa starting from the group generated by three half-turns $\Gamma = \langle a, b, c \rangle$ one can easily obtain its two generator subgroup $\Gamma_0 = \langle ab, ac \rangle$. Since $\Gamma = \Gamma_0 \cup a\Gamma_0$ ($a^2 = id$), index of $\Gamma_0$ in $\Gamma$ is equal to 2 or 1 and depends on if $a$ is contained in $\Gamma_0$ or not.

## 3. Matrix representation.

In order to obtain a matrix representation for a group generated by three half-turns one has to fix the axes of the generators. We do it by means of complex distances between the hyperbolic lines ([2], p.67-70). We define the *complex distance* between two oriented hyperbolic lines as follows: its real part equals the length of the common perpendicular to the lines and imaginary part is given by the angle between the lines, taken from the first to the second line respecting orientations.

We call three complex distances between the axes of the generators by the *parameters* of the group $\Gamma = \langle a, b, c \rangle$. Let us remark that starting from this place by group $\Gamma = \langle a, b, c \rangle$ we mean the group with the given generators taken in fixed order, so, in fact, we work with a *marked group* $\Gamma$ generated by $a, b, c$.

$$par(\langle a, b, c \rangle) = (\mu(a,b),\ \mu(a,c),\ \mu(c,b)) = (\mu_0,\ \mu_1,\ \mu_2).$$

(Here $\mu(\ldots,\ldots) = x + iy$, $x \in \mathbb{R}^+$, $y \in [0; 2\pi)$ — the complex distance between oriented axes of the half-turns.)

**Remark.** Since the directions of the half-turns are undefined the parameters are defined up to orientation of the axes. We can canonically link the orientation of a hyperbolic line with the matrix of the half-turn in $SL(2,\mathbb{C})$ (see [2], p.63). This means that fixing the directions of the axes of generators equals choosing the inverse image after canonical projection $p: SL(2,\mathbb{C}) \to PSL(2,\mathbb{C})$ of our $PSL(2,\mathbb{C})$-subgroup in $SL(2,\mathbb{C})$. Following common agreement we usually call this inverse image by representation of our group in $SL(2,\mathbb{C})$.

**Lemma 3.1.** *Complex vector $par(\langle a, b, c \rangle)$ determines group $\langle a, b, c \rangle$ uniquely up to conjugation in $PSL(2, \mathbb{C})$, i.e. if $par(\langle a, b, c \rangle) = par(\langle a', b', c' \rangle)$ then there exists $h \in PSL(2, \mathbb{C})$ such that $a' = hah^{-1}$, $b' = hbh^{-1}$, $c' = hch^{-1}$.*

PROOF: Consider the axes of the half-turns $a, b, c$ in $\mathbb{H}^3$. We endow the axes by orientations and join them in couples by common perpendiculars. The result is a hyperbolic right-angled hexagon:

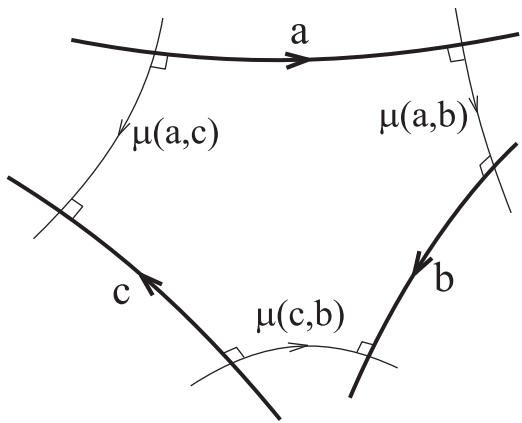

**Figure 1.**

Vector $par(\langle a, b, c \rangle)$ represents three complex lengths of the sides of our hexagon which are common perpendiculars of the given axes. By theorem from [2] (p. 94) three pairwise non-adjacent sides define right-angled hexagon uniquely up to simultaneous change of orientations of the other three sides. It means that any other hexagon defined by another group with the same parameters can be translated to the given one by a hyperbolic isometry. This isometry obviously defines the required conjugation $h \in PSL(2, \mathbb{C})$.

It may seem that we have missed the possibility of changing of the orientations of the axes but we have already remarked that the orientations effect only representation of the group in $\mathrm{SL}(2,\mathbb{C})$ but not isometries themselves. So in this situation we can change the orientations as it is needed. □

Next thing to do is to find out some convenient representation for $\langle a,b,c\rangle$ in $\mathrm{SL}(2,\mathbb{C})$ and describe it in terms of $par(\langle a,b,c\rangle)$. We denote the matrices corresponding to the half-turns in $a$, $b$, $c$ by capitals $A$, $B$, $C$.

Half-turn in oriented hyperbolic line $m$ is represented by a matrix $M \in \mathrm{SL}(2,\mathbb{C})$ of the form $\begin{pmatrix} m_{11} & m_{12} \\ m_{21} & -m_{11} \end{pmatrix}$, which is characterized by $M^2 = -I$ or equivalently $\mathrm{tr}(M) = 0$. This matrix $M$ is called *normalized line-matrix* of the line $m$. By extending transformation $M$ to the absolute $\overline{\mathbb{C}} = \partial \mathbb{H}^3$ it is easy to find the set of its fixed points in $\overline{\mathbb{C}}$ denoted by $\mathit{fix}(M)$:

$$(m_{11}z + m_{12})(m_{21}z - m_{11})^{-1} = z;$$

$$m_{21}z^2 - 2m_{11}z - m_{12} = 0.$$

The roots $\{z_1, z_2\} = \mathit{fix}(M)$ of the equation are the ends of the arc orthogonal to $\mathbb{C}$ which is the line $m$ in $\mathbb{H}^3$.

After a suitable conjugation of the group $\langle A,B,C\rangle$ in $\mathrm{SL}(2,\mathbb{C})$ one can suppose that $\mathit{fix}(A) = \{1/\beta; -1/\beta\}$, $\mathit{fix}(B) = \{\beta; -\beta\}$ ($\beta \in \mathbb{C}$ and $|\beta| \geq 1$):

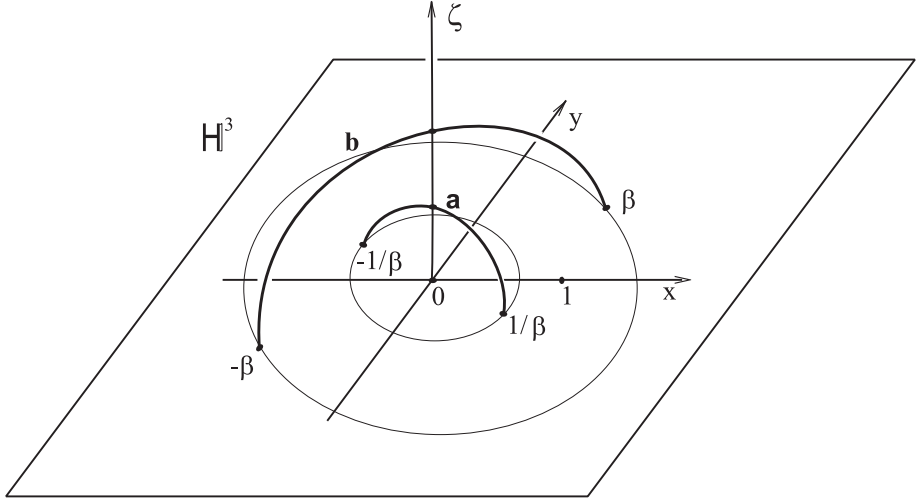

**Figure 2.**

We have:

$$A = \begin{pmatrix} 0 & i/\beta \\ i\beta & 0 \end{pmatrix}, \quad B = \begin{pmatrix} 0 & i\beta \\ i/\beta & 0 \end{pmatrix}, \quad C = \begin{pmatrix} c_{11} & c_{12} \\ c_{21} & -c_{11} \end{pmatrix},$$

where $-c_{11}^2 + c_{12}c_{21} = 1$ since $C \in \mathrm{SL}(2,\mathbb{C})$.

To find $\beta$ and $c_{ij}$ in terms of $par(\langle a,b,c\rangle)$ we use the following important formula, which allows to find the complex distance between hyperbolic lines $m_1$, $m_2$ from the matrices $M_1$, $M_2$ of the lines ([2], p.68):

$$cosh(\mu(m_1, m_2)) = -\frac{1}{2}\text{tr}(M_1 M_2).$$

Denote $\rho_k = -2\cosh\mu_k$ ($k = 0, 1, 2$). Then

$$\rho_0 = \text{tr}(AB) = \text{tr}\begin{pmatrix} -1/\beta^2 & 0 \\ 0 & -\beta^2 \end{pmatrix};$$

$$\rho_1 = \text{tr}(AC) = \text{tr}\begin{pmatrix} ic_{21}/\beta & -ic_{11}/\beta \\ ic_{11}\beta & ic_{12}\beta \end{pmatrix};$$

$$\rho_2 = \text{tr}(BC) = \text{tr}\begin{pmatrix} ic_{21}\beta & -ic_{11}\beta \\ ic_{11}/\beta & ic_{12}/\beta \end{pmatrix}.$$

We obtain system of equations on $\beta$ and $c_{ij}$:

$$\begin{cases} -1/\beta^2 - \beta^2 = \rho_0; \\ ic_{21}/\beta + ic_{12}\beta = \rho_1; \\ ic_{21}\beta + ic_{12}/\beta = \rho_2; \\ c_{11} = i\sqrt{c_{12}c_{21} + 1}. \end{cases}$$

Solving this equations we find out all parameters in our representation. Let us remark that this system has not unique solution but it really does not matter, because different solutions give representations of the groups $\Gamma = \langle a, b, c \rangle$ which are conjugate in $PSL(2, \mathbb{C})$. The last statement is a direct consequence of Lemma 3.1 since by the construction the different solutions correspond to the same $par(\Gamma)$. For definiteness we can always fix the analitic branches of the square roots in the following formulae.

Now we can write down the representation of $\langle a, b, c \rangle$ in $SL(2, \mathbb{C})$:

$$A = \begin{pmatrix} 0 & i/\beta \\ i\beta & 0 \end{pmatrix}, \quad B = \begin{pmatrix} 0 & i\beta \\ i/\beta & 0 \end{pmatrix}, \quad C = \begin{pmatrix} c_{11} & c_{12} \\ c_{21} & -c_{11} \end{pmatrix};$$

$$\beta = \sqrt{\frac{-\rho_0 + \sqrt{\rho_0^2 - 4}}{2}},$$

$$c_{21} = \frac{\rho_1/\beta - \rho_2\beta}{i/\beta^2 - i\beta^2}, \quad c_{12} = \frac{-\rho_1\beta + \rho_2/\beta}{i/\beta^2 - i\beta^2}, \quad c_{11} = i\sqrt{c_{12}c_{21} + 1};$$

$$\rho_k = -2\cosh(\mu_k) \text{ for } k = 0, 1, 2.$$

To obtain the representation of $\langle a, b, c \rangle$ in $PSL(2, \mathbb{C})$ one has to consider the image of the representation in $SL(2, \mathbb{C})$ under canonical projection $p : SL(2, \mathbb{C}) \to PSL(2, \mathbb{C})$.

## 4. Arithmeticity test.

In the following definition of arithmetic Kleinian group we consider a generalized quaternion algebra A over a number field $k \leq \mathbb{C}$. We denote by $A^1$ the group of quaternions of norm 1. The algebra $A \otimes \mathbb{C}$ is isomorphic to $M(2, \mathbb{C})$ and we shall think of $A$ as a $k$-subalgebra of $\mathbb{C}$-algebra $M(2, \mathbb{C})$. The reduced trace and norm of $A$ coincide with the usual trace and determinant functions in $M(2, \mathbb{C})$.

**Definition 4.1.** ([8], [4]) *A subgroup $\Gamma$ of* $\mathrm{PSL}(2,\mathbb{C})$ *is called* arithmetic *if there exist*

*(1) a quaternion algebra $A$ over a number field $k \leq \mathbb{C}$ which has exactly one complex place (pair of complex conjugate embeddings), and such that $A$ is* ramified *at the real places: i.e. $A \otimes_\phi \mathbb{R} = \mathbb{H}$ (the Hamilton quaternions) for every real place $\phi : k \to \mathbb{R}$.*

*(2) an order $O$ of $A$ such that $\Gamma$ is* commensurable *with $p(O^1)$, where $O^1 = O \cap A^1$ and $p : \mathrm{SL}(2,\mathbb{C}) \to \mathrm{PSL}(2,\mathbb{C})$ — canonical projection (i.e. there exists a subgroup $\Gamma_1$ of $\Gamma$ of finite index such that $\Gamma_1$ is a subgroup of $p(O^1)$ of finite index).*

**Definition 4.2.** [10] *A subgroup $\Gamma$ of* $\mathrm{PSL}(2,\mathbb{C})$ *is called* nearly arithmetic *if*
*(1) $\Gamma$ is a subgroup of an arithmetic Kleinian group.*
*(2) $\Gamma$ is not isomorphic to a free product of cyclic groups.*

Let $\Gamma$ be a non-elementary finitely generated subgroup of $\mathrm{SL}(2,\mathbb{C})$ and let $\mathbb{Q}(\mathrm{tr}\Gamma) = \mathbb{Q}(\{\mathrm{tr}\ g : g \in \Gamma\})$ denote the trace field of $\Gamma$. Then

$$A\Gamma = \{\sum a_i g_i : a_i \in \mathbb{Q}(\mathrm{tr}\Gamma),\ g_i \in \Gamma\}$$

is a quaternion algebra over $\mathbb{Q}(\mathrm{tr}\Gamma)$ [8]. Futhermore, if we set $\Gamma^{(2)} = \langle g^2 \mid g \in \Gamma \rangle$, then

$$k\Gamma = \mathbb{Q}(\mathrm{tr}\Gamma^{(2)}) = \mathbb{Q}(\{\mathrm{tr}^2 g : g \in \Gamma\})$$

is an invariant of the commensurability class of $\Gamma$ as is the quaternion algebra $A\Gamma^{(2)}$ [8]. For a group $\Gamma \leq \mathrm{PSL}(2,\mathbb{C})$ we use the same notation but strictly we are referring to the inverse image of $\Gamma$ in $\mathrm{SL}(2,\mathbb{C})$.

If $\Gamma$ is arithmetic Kleinian group, then the invariant trace field $k\Gamma$ and the quaternion algebra $A\Gamma^{(2)}$ recapture the arithmetic structure defining $\Gamma$ so that arithmetic Kleinian groups can be characterised intrinsically. More precisely (see [8], [4]):

**Theorem 4.1.** *A finitely generated subgroup $\Gamma$ of $PSL(2,\mathbb{C})$ is nearly arithmetic if and only if the following conditions are satisfied:*
*(1) $\mathrm{tr}\Gamma$ consists of algebraic integers;*
*(2) $k\Gamma$ is a number field which has exactly one complex place;*
*(3) $A\Gamma^{(2)}$ is ramified at all real places;*
*(4) $\Gamma$ is not isomorphic to a free product of cyclic groups.*
*Furthermore, a nearly arithmetic group is arithmetic if and only if it has a finite covolume.*

Now let us return to the groups generated by three half-turns. After some common calculations with traces ([8], [4]) it is easy to see that if we denote $g_0 = ab$, $g_1 = ac$ then for $\Gamma = \langle a, b, c \rangle$ and $\Gamma' = \langle g_0, g_1 \rangle$

$$\mathrm{tr}\Gamma' = \mathbb{Q}(\mathrm{tr}(g_0),\ \mathrm{tr}(g_1),\ \mathrm{tr}(g_0 g_1)) = \mathbb{Q}(\rho_0,\ \rho_1,\ \rho_2);$$

$$k\Gamma = k\Gamma' = \mathbb{Q}(\mathrm{tr}(g_0^2),\ \mathrm{tr}(g_1^2),\ \mathrm{tr}(g_0^2 g_1^2)) = \mathbb{Q}(\rho_0^2,\ \rho_1^2,\ \rho_0\rho_1\rho_2);$$

$$A\Gamma^{(2)} = \left(\frac{\mathrm{tr}(g_0^2)^2 - 4, \mathrm{tr}[g_0^2, g_1^2] - 2}{k\Gamma}\right) = \left(\frac{\rho_0^2(\rho_0^2 - 4), \rho_0^2 \rho_1^2(\rho_0^2 + \rho_1^2 + \rho_2^2 - \rho_0\rho_1\rho_2 - 4)}{k\Gamma}\right).$$

So we can reformulate the previous theorem for the groups generated by three half-turns:

**Theorem 4.2. (Arithmeticity test)** *Subgroup $\Gamma = \langle a, b, c \rangle$ of $PSL(2, \mathbb{C})$ generated by three half-turns is nearly arithmetic if and only if:*

*(1) $\rho_0$, $\rho_1$, $\rho_2$ are algebraic integers;*
*(2) $k\Gamma = \mathbb{Q}(\rho_0^2, \rho_1^2, \rho_0\rho_1\rho_2)$ has exactly one complex place;*
*(3) for every real place $\phi$ of $k\Gamma$:*
$$\phi(\rho_0^2(\rho_0^2 - 4)) < 0, \ \phi(\rho_0^2\rho_1^2(\rho_0^2 + \rho_1^2 + \rho_2^2 - \rho_0\rho_1\rho_2 - 4)) < 0;$$
*(4) $\Gamma$ is not isomorphic to a free product of cyclic groups.*
*$\Gamma$ is arithmetic if and only if it is nearly arithmetic and has a finite covolume.*

## 5. Regular triangle groups.

In this section we discuss a partial case of the groups generated by three half-turns. We assume that all the three complex distances between the axes of the generators are equal and call corresponding groups *regular triangle groups*. For this class of groups we prove that it contains only finitely many arithmetic groups and present a list containing all non cocompact arithmetic regular triangle groups.

For the regular triangle groups matrix representation from section 3 can be rewritten in simpler form. Here $par(\langle a, b, c \rangle) = (\mu, \mu, \mu)$ and we have:

$$A = \begin{pmatrix} 0 & i/\beta \\ i\beta & 0 \end{pmatrix}, \ B = \begin{pmatrix} 0 & i\beta \\ i/\beta & 0 \end{pmatrix}, \ C = \begin{pmatrix} c_{11} & c_{21} \\ c_{21} & -c_{11} \end{pmatrix};$$

$$\rho = -2\cosh(\mu),$$
$$\beta = \sqrt{\frac{-\rho + \sqrt{\rho^2 - 4}}{2}},$$
$$c_{21} = i\frac{1/\beta^2 + \beta^2}{1/\beta + \beta},$$
$$c_{11} = i\sqrt{c_{21}^2 + 1}.$$

We start with a technical lemma, which is essential part of proof of the finiteness theorem:

**Lemma 5.1.** *Let $\Gamma$ be a regular triangle group. Then if $|\rho| \geq \rho_* = 6.4$ then $\Gamma$ splits into free product of cyclic groups.*

PROOF: Let $C_0$, $C_1$, $C_2$ be invariant circles of $A$, $B$ and $C$ on $\bar{\mathbb{C}}$ considered as a limit plane of $\mathbb{H}^3$. When this three circles do not intersect each other, group $\Gamma$ splits into a free product of cyclics by Klein combination theorem (see [11], p. 139).

Each invariant circle of a hyperbolic isometry passes through its fixed points and we can actually assume that the points lie on the diameter of the corresponding circle in $\mathbb{C}$. Let us find the conditions on the parameters which can guarantee that the invariant circles don't intersect. We have:

$$\forall z \in C_0: \ |z| \leq 1/|\beta|,$$
$$\forall z \in C_1: \ |z| \geq |\beta|,$$
$$\exists R_1, R_2 \ | \ \forall z \in C_2: \ R_1 \leq |z| \leq R_2.$$

Radius of $C_2$:
$$r = \frac{1}{|c_{21}|} \leq \frac{1 + 1/|\beta|^2}{|\beta|(1 - 1/|\beta|^4)} \overset{\text{def}}{=} \phi(|\beta|),$$

its center in $\mathbb{C}$ is $c_{11}/c_{21} = i\sqrt{1+1/r^2}$. Hence we can take:
$$R_1 = R_1(|\beta|) = \sqrt{|1-\phi^2|} - \phi,$$
$$R_2 = R_2(|\beta|) = \sqrt{|1+\phi^2|} + \phi.$$

We obtain two functions $R_1(x)$ and $R_2(x)$ of real variable bounding the third circle. Calculating derivatives it is easy to see that when $x$ is large enough ($x \geq 2$), $R_1(x)$ monotonically increase and $R_2(x)$ monotonically decrease. It follows that if the circles do not intersect for some value $\beta_*$ of $|\beta|$ ($\beta_* \geq 2$) then they do not intersect for all $\beta$ with $|\beta| \geq \beta_*$ (see figure 3). The value of $\beta_*$ is defined by one of the equations: $R_1(x) = 1/x$ or $R_2(x) = x$ and so approximately $\beta_* = 2.49455$.

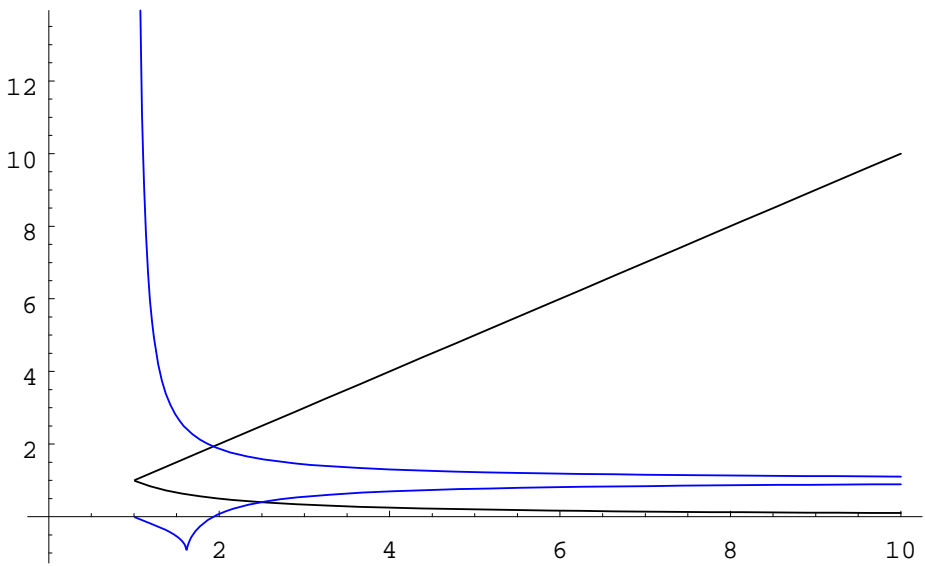

**Figure 3.**

It remains to recall that $|\rho| = |-1/\beta^2 - \beta^2| \leq 1/|\beta|^2 + |\beta|^2$ so $|\rho| \geq 1/\beta_*^2 + \beta_*^2$ implies $1/|\beta|^2 + |\beta|^2 \geq 1/\beta_*^2 + \beta_*^2$. Since function $1/x^2 + x^2$ monotonically increase when $x \geq 1$, $|\rho| \geq 1/\beta_*^2 + \beta_*^2$ follows $|\beta| \geq \beta_*$ which means that our group splits into free product. So we can take $\rho_* = 1/\beta_*^2 + \beta_*^2 = 6.4$.  □

Now we are ready to prove the finiteness theorem for the number of arithmetic regular triangle groups.

**Theorem 5.2.** *There are only finitely many conjugacy classes of nearly arithmetic regular triangle Kleinian groups.*

PROOF: By arithmeticity test (Theorem 4.2) in order group $\Gamma$ to be arithmetic $\rho = -1/\beta^2 - \beta^2$ has to be algebraic integer. Let $D(\rho)$ be discriminant of $\rho$ over $\mathbb{Q}$. From definition of discriminant it follows that
$$D = (\rho - \bar\rho)^2 \prod_{i=1}^{n}(\rho - \rho_i)^2(\bar\rho - \rho_i)^2 \prod_{1 \leq i < j \leq n}(\rho_i - \rho_j)^2$$

(here $(n+2)$ is the order of algebraic integer $\rho$ and $\rho_i$, $i=1,\ldots,n$ denote other roots of the minimum polynomial for $\rho$).

Since $\rho$ is an algebraic integer its reduced norm $N(D)$ is a non-zero rational integer (see [13]) and so
$$|N(D)| = |D| \geq 1.$$
From the other hand
$$|N(D)| \leq |\rho-\bar{\rho}|^2 \prod_{i=1}^{n} |(\rho-\rho_i)(\bar{\rho}-\rho_i)|^2 \prod_{1\leq i<j\leq n}(\rho_i-\rho_j)^2 \leq 4\rho_*^2 K^{2n} \prod_{1\leq i<j\leq n}(\rho_i-\rho_j)^2,$$
$$K = |\rho_*|^2 + 4|\rho_*| + 4 = 70.56.$$
(We applied Lemma 5.1 and used conditions (2) and (3) of the arithmeticity test from section 4, stating that if $\Gamma$ is nearly arithmetic then $\rho_i$ are real numbers and $|\rho_i| \leq 2$, $i=1,\ldots,n$.)

To continue this chain of inequalities we recall one number theoretical result (see [10] for discussion):

If $-1 \leq x_1 < x_2 < \cdots < x_r \leq 1$ where $r \geq 3$ then
$$\prod_{1\leq i<j\leq r}(x_i-x_j)^2 \leq M_r = \frac{2^2 3^3 \cdots r^r 2^2 3^3 \cdots (r-2)^{r-2}}{3^3 5^5 \cdots (2r-3)^{2r-3}}$$
and $M_n^{1/n(n-1)} \to 1/4$ while $n \to \infty$. \hfill (*)

In our case it gives
$$N(D) \leq 4\rho_*^2 K^{2n} 2^{n(n-1)} M_n$$
and (*) implies $N(D) \leq 4\rho_*^2 K^{2n} 2^{n(n-1)} M_n \sim K^{2n} 2^{n^2} 2^{-n^2} \to 0$, $n \to \infty$. Hence there exists number $n_0$ such that $\forall n \geq n_0$: $N(D) < 1$ and so $\rho$ can not be algebraic integer if its order is greater then $n_0$. By Lemma 5.1 if the group do not split into free product then $|\rho| \leq \rho_*$. It remains to recall a basic number-theoretic fact that there are only finitely many algebraic integers of bounded degree such that all their conjugates are bounded in absolute value. $\square$

**Corollary 5.3.** *There are only finitely many conjugacy classes of arithmetic regular triangle Kleinian groups.*

**Remark.** Actually, all the regular triangle groups are commensurable with the groups generated by an elliptic rotation of order 3 and a half-turn. It means that Theorem 5.2 can be considered as a corollary from C. Maclachlan and G.J. Martin result [10] on the finiteness of the number of arithmetic Kleinian groups generated by two elliptic elements. The main idea of the Theorem 5.2 is to demonstrate how does our construction work in the simpliest possible situation. Now the same method can be applied to the other, more complicated cases and, finally, it gives a new approach to the general problem of classifying all two-generator arithmetic Kleinian groups.

We finish this section with more careful discussion of the non cocompact regular triangle groups. It is well known that every non cocompact arithmetic group is commensurable with some Bianchi group [3]. This means that in non cocompact case the field of definition $k\Gamma$ of an arithmetic group $\Gamma$ is $\mathbb{Q}(\sqrt{-d})$ ($d$ is a positive integer). It follows that parameter $\rho$ is an algebraic integer of order 2, so $\rho = k/2 + \sqrt{-d}/2$, where $k \in \mathbb{Z}$, $d \in \mathbb{Z}^+$ and $k^2 + d$ is devisible by 4. In this case we can apply the estimate from Lemma 5.1 to write down the list of all possible values for $\rho$ giving nearly arithmetic $\Gamma$. After more precise computer calculations with invariant circles we have reduced this list to 48 representatives (Table 1). For further refinements one can use the specific symmetries of the regular triangle groups from one side and additional information about Bianchi groups and corresponding orbifolds (see [6], [7], [1]) from the other. Finally we have obtained three examples (see Figures 4,5,6) of non cocompact arithmetic regular triangle groups and our conjecture was that there are no other groups in this class. Lately Colin Maclachlan told me that this conjecture follows from his joint result with Gaven Martin on non cocompact arithmetic Kleinian groups generated by two elliptic elements which will be published soon. Concerning this, let us remark once more that the connection of the generalized triangle groups with the Kleinian groups generated by two elliptics, which takes place in the regular case, does not hold in general situation.

On the following pictures we present the singular sets of the orbifolds $\mathbb{H}^3/\Gamma$ where $\Gamma = \langle a, b, c \rangle$ is an arithmetic non cocompact regular triangle group. The singular sets are knoted graphs with singular values denoted on the pictures (the value is 2 if it is not denoted). The underlying spaces of the orbifolds are punctured 3-spheres $\mathbb{S}^3$.

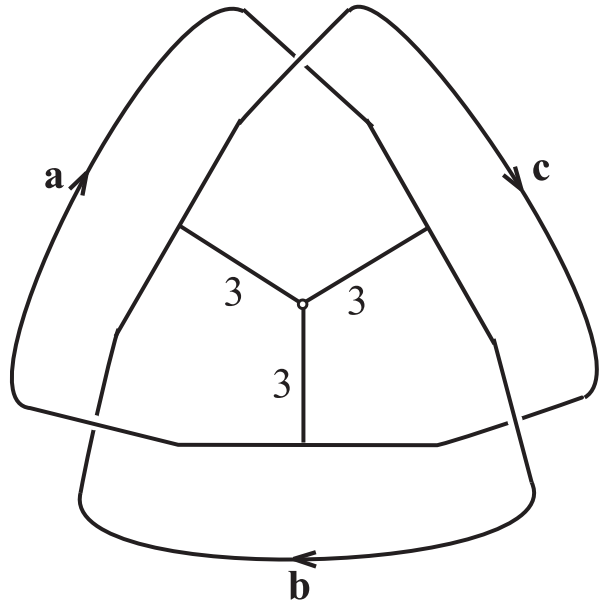

$\Gamma = \langle a, b, c \mid a^2, b^2, c^2, w(a, b, c), w(b, c, a), w(c, a, b) \rangle$, $w(x, y, z) = z^{-1}y^{-1}x^{-1}yzx$;
$\rho = -3/2 + \sqrt{-3}/2$, $k\Gamma = \mathbb{Q}(\sqrt{-3})$.

**Figure 4.**

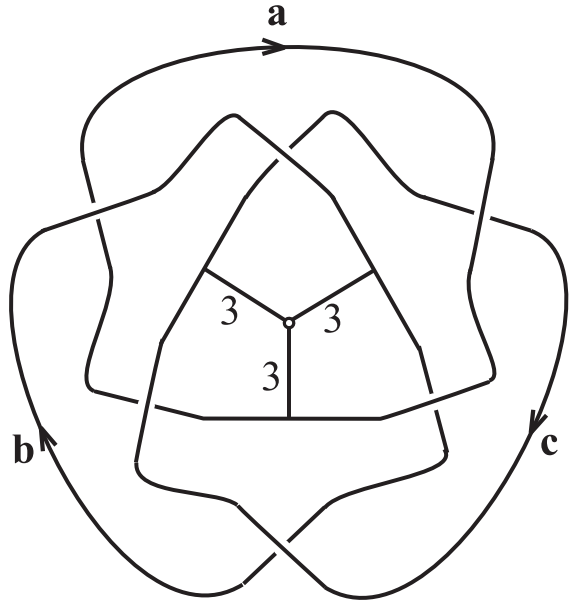

$\Gamma = \langle a, b, c \mid a^2, b^2, c^2, w(a,b,c), w(b,c,a), w(c,a,b) \rangle, \ w(x,y,z) = x^{-1}y^{-1}zy;$
$\rho = -1/2 + \sqrt{-3}/2, \ k\Gamma = \mathbb{Q}(\sqrt{-3}).$
**Figure 5.**

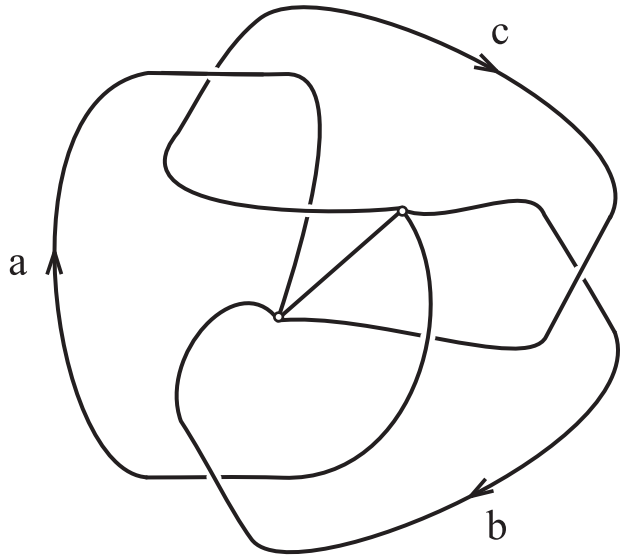

$\Gamma = \langle a, b, c \mid a^2, b^2, c^2, w(a,b,c), w(b,c,a), w(c,a,b) \rangle,$
$w(x,y,z) = y^{-1}x^{-1}yz^{-1}y^{-1}zx^{-1}z^{-1}x;$
$\rho = 1/2 + \sqrt{-3}/2, \ k\Gamma = \mathbb{Q}(\sqrt{-3}).$
**Figure 6.**

## 6. General case.

In this section we will discuss how our methods could be applied to the general problem of classifying all arithmetic Kleinian group generated by three half-turns and what are the difficulties.

The idea is first to prove some analog of Lemma 5.1 and then apply the arguments of the proof of Theorem 5.2 to each of the three parameters.

**Lemma 6.1. (Conjecture)** *Let $\Gamma$ be a Kleinian group generated by three half-turns. Then when $|\rho_i| \geq \rho_* = 6.4$ $(i = 0, 1, 2)$ $\Gamma$ splits into free product of cyclic groups.*

This Lemma has appeared as a result of numerical experiments with invariant circles. We can prove it in partial cases, but we are not ready yet to give the complete analitic proof.

The main problem is that Lemma 6.1 is not enough to prove the finiteness theorem. In order to do it the same way as earlier we need more powerfull result, that there is a ball in $\mathbb{C}^3$ containing all possible values of the parameters $\rho_i$ for which the group does not split into free product, but this seems not to be true. There arise some interesting exceptional sequences of parameters, for example:

$$\{\rho_0 = n/2 + \sqrt{-3}/2,\ \rho_1 = \rho_2 = 1/2 + \sqrt{-3}/2\}_{n \in \mathbb{Z}}$$

For all the groups in this sequence the invariant circles of $A$, $B$ and $C$ intersect, so we can not use the Klein combination theorem. From one hand it could be some other way to avoid this sort of sequences, for example, using the finiteness of volume arguments. From the other hand it might be interesting to investigate the groups in this sequence and their factor orbifolds (all this groups are certainly discrete because they are nearly arithmetic or split into free products).

Sobolev Institute of Mathematics
Novosibirsk, 630090, Russia
email: mbel@math.nsc.ru

| N | $\rho$ | Field | N | $\rho$ | Field |
|---|---|---|---|---|---|
| 1 | $-3+\sqrt{-1}$ | $\mathbb{Q}(\sqrt{-1})$ | 25 | $-1+\sqrt{-5}$ | $\mathbb{Q}(\sqrt{-5})$ |
| 2 | $-2+\sqrt{-1}$ | $\mathbb{Q}(\sqrt{-1})$ | 26 | $\sqrt{-5}$ | $\mathbb{Q}(\sqrt{-5})$ |
| 3 | $-1+\sqrt{-1}$ | $\mathbb{Q}(\sqrt{-1})$ | 27 | $-1+\sqrt{-6}$ | $\mathbb{Q}(\sqrt{-6})$ |
| 4 | $\sqrt{-1}$ | $\mathbb{Q}(\sqrt{-1})$ | 28 | $\sqrt{-6}$ | $\mathbb{Q}(\sqrt{-6})$ |
| 5 | $1+\sqrt{-1}$ | $\mathbb{Q}(\sqrt{-1})$ | 29 | $-5/2+\sqrt{-7}/2$ | $\mathbb{Q}(\sqrt{-7})$ |
| 6 | $2+\sqrt{-1}$ | $\mathbb{Q}(\sqrt{-1})$ | 30 | $-3/2+\sqrt{-7}/2$ | $\mathbb{Q}(\sqrt{-7})$ |
| 7 | $-2+\sqrt{-4}$ | $\mathbb{Q}(\sqrt{-1})$ | 31 | $-1/2+\sqrt{-7}/2$ | $\mathbb{Q}(\sqrt{-7})$ |
| 8 | $-1+\sqrt{-4}$ | $\mathbb{Q}(\sqrt{-1})$ | 32 | $1/2+\sqrt{-7}/2$ | $\mathbb{Q}(\sqrt{-7})$ |
| 9 | $\sqrt{-4}$ | $\mathbb{Q}(\sqrt{-1})$ | 33 | $3/2+\sqrt{-7}/2$ | $\mathbb{Q}(\sqrt{-7})$ |
| 10 | $2+\sqrt{-4}$ | $\mathbb{Q}(\sqrt{-1})$ | 34 | $-1+\sqrt{-7}$ | $\mathbb{Q}(\sqrt{-7})$ |
| 11 | $-2+\sqrt{-2}$ | $\mathbb{Q}(\sqrt{-2})$ | 35 | $-5/2+\sqrt{-11}/2$ | $\mathbb{Q}(\sqrt{-11})$ |
| 12 | $-1+\sqrt{-2}$ | $\mathbb{Q}(\sqrt{-2})$ | 36 | $-3/2+\sqrt{-11}/2$ | $\mathbb{Q}(\sqrt{-11})$ |
| 13 | $\sqrt{-2}$ | $\mathbb{Q}(\sqrt{-2})$ | 37 | $-1/2+\sqrt{-11}/2$ | $\mathbb{Q}(\sqrt{-11})$ |
| 14 | $2+\sqrt{-2}$ | $\mathbb{Q}(\sqrt{-2})$ | 38 | $1/2+\sqrt{-11}/2$ | $\mathbb{Q}(\sqrt{-11})$ |
| 15 | $-5/2+\sqrt{-3}/2$ | $\mathbb{Q}(\sqrt{-3})$ | 39 | $3/2+\sqrt{-11}/2$ | $\mathbb{Q}(\sqrt{-11})$ |
| 16 | $-3/2+\sqrt{-3}/2$ | $\mathbb{Q}(\sqrt{-3})$ | 40 | $-3/2+\sqrt{-15}/2$ | $\mathbb{Q}(\sqrt{-15})$ |
| 17 | $-1/2+\sqrt{-3}/2$ | $\mathbb{Q}(\sqrt{-3})$ | 41 | $-1/2+\sqrt{-15}/2$ | $\mathbb{Q}(\sqrt{-15})$ |
| 18 | $1/2+\sqrt{-3}/2$ | $\mathbb{Q}(\sqrt{-3})$ | 42 | $1/2+\sqrt{-15}/2$ | $\mathbb{Q}(\sqrt{-15})$ |
| 19 | $3/2+\sqrt{-3}/2$ | $\mathbb{Q}(\sqrt{-3})$ | 43 | $-3/2+\sqrt{-19}/2$ | $\mathbb{Q}(\sqrt{-19})$ |
| 20 | $-2+\sqrt{-3}$ | $\mathbb{Q}(\sqrt{-3})$ | 44 | $-1/2+\sqrt{-19}/2$ | $\mathbb{Q}(\sqrt{-19})$ |
| 21 | $-1+\sqrt{-3}$ | $\mathbb{Q}(\sqrt{-3})$ | 45 | $1/2+\sqrt{-19}/2$ | $\mathbb{Q}(\sqrt{-19})$ |
| 22 | $\sqrt{-3}$ | $\mathbb{Q}(\sqrt{-3})$ | 46 | $-3/2+\sqrt{-23}/2$ | $\mathbb{Q}(\sqrt{-23})$ |
| 23 | $1+\sqrt{-3}$ | $\mathbb{Q}(\sqrt{-3})$ | 47 | $-1/2+\sqrt{-23}/2$ | $\mathbb{Q}(\sqrt{-23})$ |
| 24 | $-1/2+\sqrt{-27}/2$ | $\mathbb{Q}(\sqrt{-3})$ | 48 | $1/2+\sqrt{-23}/2$ | $\mathbb{Q}(\sqrt{-23})$ |

**Table 1:** The list containing all nearly arithmetic regular triangle groups defined over $\mathbb{Q}(\sqrt{-d})$.